\documentclass[12pt,psamsfonts]{amsart}
\usepackage{amsmath}
\usepackage{amsthm}
\usepackage{amssymb}
\usepackage{amscd}
\usepackage{amsfonts}
\usepackage{amsbsy}

\newtheorem {theorem} {Theorem}

\begin{document}

\title[A counterexample to the composition condition]
{A counterexample to the composition condition conjecture for polynomial Abel differential equations}

\author[Jaume Gin\'e, Maite Grau and Xavier Santallusia]
{Jaume Gin\'e, Maite Grau and Xavier Santallusia}

\address{Departament de Matem\`atica, Inspires Research Centre, Universitat de Lleida,
Avda. Jaume II, 69; 25001 Lleida, Catalonia, Spain}

\email{gine@matematica.udl.cat}
\email{mtgrau@matematica.udl.cat}
\email{xivaer@gmail.com}

\thanks{The authors are partially supported by a
MINECO/FEDER grant number MTM2014-53703-P and by an AGAUR (Generalitat de
Catalunya) grant number 2014SGR 1204.}

\subjclass[2010]{Primary 34C25. Secondary 34C07.}

\keywords{Abel equations, center problem, composition condition, moment conditions, composition conjecture}
\date{}
\dedicatory{}

\maketitle

\begin{abstract}
The Polynomial Abel differential equations are considered a model problem for the classical Poincar\'e center--focus problem for planar polynomial systems of ordinary differential equations. Last decades several works pointed out that all the centers of the polynomial Abel differential equations satisfied the composition conditions (also called universal centers). In this work we provide a simple counterexample to this conjecture.
\end{abstract}

\section{Introduction \label{sect1}}

These last decades some authors consider polynomial Abel differential equations as a model to tackle the center problem for a trigonometric Abel differential equation coming from a planar polynomial systems of ordinary differential equations, see \cite{BFY1,BFY2,BFY3}. We denote as a {\em polynomial Abel differential equation} an ordinary differential equation of the form
\begin{equation}\label{eq3}
\frac{d y}{dx}\, =\, p(x) y^2 + q(x) y^3,
\end{equation}
where $y$ is real, $x$ is a real independent variable considered in a real interval $[a,b]$ and $p(x)$ and $q(x)$ are real polynomials in $\mathbb{R}[x]$. The {\em center problem for a polynomial Abel equation} (\ref{eq3}) is to characterize when all the solutions in a neighborhood of the solution $y=0$ take the same value when $x=a$ and $x=b$, i.e. $y(a)=y(b)$. In this framework, given any real continuous function $c(x)$, we denote by $\tilde{c}(x) \, := \, \int_{a}^{x} c(\sigma) d \sigma$ and we will say that a real continuous function $w(x)$ is {\em periodic in $[a,b]$} if $w(a)=w(b)$.
\par
Alwash and Lloyd in \cite{AL} provided a sufficient condition for an Abel trigonometric equation
\begin{equation}\label{eq2}
\frac{d \rho}{d \theta}\,  =\, a_1(\theta) \rho^2 + a_2 (\theta) \rho^3,
\end{equation}
where $\rho$ is real, $\theta$ is a real and periodic independent variable with $\theta \in [0,2\pi]$, and $a_1(\theta)$ and $a_2(\theta)$ are real trigonometric polynomials, to have a center in $[0,2\pi]$. We recall that the center problem for equation (\ref{eq2}) is to characterize when all the solutions in a neighborhood of the solution $\rho=0$ are periodic of period $2\pi$. Inspired by this work, Briskin, Fran\c{c}oise and Yomdin in \cite{BFY1} provided the following sufficient condition for the polynomial Abel equation (\ref{eq3}).
\begin{theorem}{\cite{BFY1}}\label{t1}
If there exists a real differentiable function $w$ periodic in $[a,b]$ and such that
\[ \tilde p(x)\, =\, p_1(w(x)) \quad \mbox{and} \quad \tilde q(x)\, =\, q_1(w(x))\] for some real differentiable functions $p_1$ and $q_1$,
then the polynomial Abel equation {\rm (\ref{eq3})} has a center in $[a,b]$.
\end{theorem}
In \cite{CGM3} it is shown that if the sufficient condition stated in Theorem \ref{t1} is satisfied then there is a countable set of definite integrals which need to vanish. In \cite{CGM3} it is also shown that this is equivalent to the existence of a real polynomial $w(x)$ with $w(a)=w(b)$ and two real polynomials $p_1(x)$ and $q_1(x)$ such that $\tilde{p}(x) \, = \, p_1(w(x))$ and $\tilde{q}(x) \, = \, q_1(w(x))$. This sufficient condition is known as the {\em composition condition}. \par To see that the composition condition implies that equation (\ref{eq3}) has a center in $[a,b]$ one can consider the transformation
$y(x)=Y(w(x))$ in equation (\ref{eq3}) in order to obtain the following Abel differential equation
\begin{equation}\label{eq4}
\frac{d Y}{dw} \, =\,  p_1'(w) Y^2 + q_1'(w) Y^3.
\end{equation}
Hence, there is a bijection between the solutions $Y=Y(w)$ of equation (\ref{eq4}) and the solutions $y=Y(w(x))$ of equation (\ref{eq3}). Since $w$ is periodic in $[a,b]$, we get that equation (\ref{eq3}) has a center in $[a,b]$ because $y(a)=Y(w(a))=Y(w(b))=y(b)$.

\smallskip

It turns out that all the known polynomial Abel differential equations which have a center in $[a,b]$ satisfy the composition condition. Hence in several works was established what is know as {\em composition conjecture}, see \cite{A0,GGX2} and references therein. This conjecture says that the sufficient condition given in Theorem \ref{t1} is also necessary. That is, if a polynomial Abel equation (\ref{eq3}) has a center in $[a,b]$, the conjecture states that the composition condition is satisfied. In fact in \cite{GGX} this conjecture was proved for lower degrees of the polynomials $p(x)$ and $q(x)$ of equation (\ref{eq3}). Moreover is satisfied under certain restrictions of the coefficients of the polynomial Abel differential equation, see for instance \cite{BRY}, Theorem 2 in \cite{A2}, Theorem 2 in \cite{BlY} and Theorem 7 in \cite{GGX}.

\smallskip

For a trigonometric Abel differential equation (\ref{eq2}), Alwash in \cite{A} showed that this conjecture is not true, see also \cite{A2,CGM,CGM2,GGX0}. The {\em composition condition} for a trigonometric Abel differential equation (\ref{eq2}) is that there exist real polynomials $p_1(x), p_2(x) \in \mathbb{R}[x]$ and a trigonometric polynomial $\omega(\theta)$ such that $\tilde{a}_i(\theta) \, = \, p_i(\omega(\theta))$, for $i=1,2$. Recall that $\tilde{a_i}(\theta) \, := \, \int_{0}^{\theta} a_i(s) ds$. The fact that $\omega(\theta)$ and $p_1, p_2$ can be taken to be polynomials is proved in \cite{CGM3,GGL}. There exist several counterexamples of the fact that the composition conjecture is not satisfied in the trigonometric case. The authors of \cite{A,A2,CGM} provide examples of trigonometric polynomials $a_1(\theta)$ and $a_2(\theta)$ for which the corresponding trigonometric Abel differential equation (\ref{eq2}) has a center and does not satisfy the composition condition. A survey of the last results for polynomial and trigonometric Abel equations is given in \cite{GGX2}.

\smallskip

The main result of this note is the following.

\begin{theorem}\label{th1}
The polynomial Abel equation (\ref{eq3}) with
\[
p(x)=40 x^4 - 30 x^2 + 2, \quad q(x)=75 x^9 - 150 x^7 + 88 x^5 - 10 x^3 - 3 x,
\]
has a center and does not satisfy the composition condition.
\end{theorem}

In the following section we proof the main result of this note.

\section{Proof of Theorem \ref{th1} \label{sect2}}

System (\ref{eq3}) with $p(x)$ and $q(x)$ given by statement of the theorem has the invariant algebraic curves
\[
\begin{array}{ll}
f_1=& (2 + 8 x y - 24 x^3 y + 16 x^5 y + y^2 + 2 x^2 y^2 - 34 x^4 y^2 \\
   & + 88 x^6 y^2 - 87 x^8 y^2 + 30 x^{10} y^2)/2,\\
f_2=& (3 + 12 x y - 42 x^3 y + 30 x^5 y + 2 y^2 + 3 x^2 y^2 - 72 x^4 y^2 +\\
   & 202 x^6 y^2 - 210 x^8 y^2 + 75 x^{10} y^2)/3,\\
f_3=&1 + 3 x y - 8 x^3 y + 5 x^5 y,
\end{array}
\]
and the rational first integral $H(x,y)=y^2f_1^3/f_2^4$. This first integral satisfies that
$H(1,y)=H(-1,y)=y^2$ consequently this Abel trigonometric equation has a center. Moreover attending to the first integral obtained system (\ref{eq3}) with $p(x)$ and $q(x)$ given by statement of the theorem admits a type of first integral studied in \cite{GS} for Abel equations. In particular corresponds to a case with five solutions, that is, $n=5$. In order to prove that this Abel equation does not satisfies the composition condition we must to recall the equivalence between composition condition and the existence of a universal center, see \cite{GGL}.

An explicit expression for the first return map of the differential
equation (\ref{eq3}) was given in \cite{Bru}, see also \cite{Br}.
This expression is given in terms of the following iterated
integrals, of order $k$,
\[
I_{i_1,\ldots,i_k}(\lambda):=\int \cdots \int_{0 \le s_1 \le \cdots \le
s_k \le 2 \pi} a_{i_k}(s_k)\cdots a_{i_1}(s_1)\, d s_k \cdots ds_1,
\]
where, by convention, for $k=0$ we assume that this equals $1$. Let
$y(x;y_0;\lambda)$, $x \in  [a, b]$, be the Lipschitz
solution of the differential equation (\ref{eq3}) corresponding to
a sequence $\lambda=(\lambda_1,\lambda_2,\ldots)$ of parameters of equation
(\ref{eq3}) with initial value $y(a;y_0;\lambda)=y_0$. Then
$P(\lambda)(y_0):=y(b;y_0;\lambda)$ is the first return map of this
differential equation, and in \cite{Bru, Br} it is proved the
following:

\begin{theorem}
For sufficiently small initial values $y_0$ the first return map
$P(a)$ is an absolute convergent power series $P(\lambda)(y_0)= y_0
+ \sum_{n=1}^{\infty} c_n(\lambda) y_0^{n+1}$, where
\[
c_n(\lambda)= \sum_{i_1+\cdots+i_k=n} c_{i_1, \ldots, i_k} I_{i_1, \ldots,
i_k} (\lambda), \quad \mbox{and}
\]
\[
c_{i_1, \ldots, i_k}
=(n-i_1+1)\cdot(n-i_1-i_2+1)\cdot(n-i_1-i_2-i_3+1)\cdots 1.
\]
\end{theorem}

The following definition is given in \cite{Br2}. Equation
(\ref{eq3}) determines a {\it universal center} if for all
positive integers $i_1,\ldots,i_k$ with $k \ge 1$ the iterated
integral $I_{i_1,\ldots,i_k}(\lambda)=0$. Moreover in \cite{GGL} it was
proved that equation (\ref{eq3}) has a universal center if and only
if the composition condition is satisfied.

System (\ref{eq3}) with $p(x)$ and $q(x)$ given by statement of the theorem has the iterated integral
\[
I_{122}=\int_{-1}^1 \int_{-1}^{x_3} \int_{-1}^{x_2} p(x_1) q(x_2) q(x_3) dx_1 dx_2 dx_3 = -\frac{131072}{6235515}.
\]
Hence, we have a non universal center and this completes the proof.\\

In fact there is a straightforward way to see that system (\ref{eq3}) with $p(x)$ and $q(x)$ given by statement of the theorem does not satisfies
the composition condition. This consists in to see that the integral
\[
\int_{-1}^1  \tilde p(x) \tilde q^2(x) p(x) dx
\]
is not null, where recall that now $\tilde c (x):= \, \int_{-1}^{x} c(\sigma) d \sigma$.


\begin{thebibliography}{99}

\bibitem{A} {\sc M.A.M. Alwash}, {\it On a condition for a center of cubic non-autonomous equations},  Proc. Roy. Soc. Edinburgh Sect. A  {\bf 113} (1989), 289--291.

\bibitem{A0} {\sc M.A.M. Alwash}, {\it On the composition conjectures},  Electron. J. Differential Equations {\bf 2003}, No. 69, 4 pp.

\bibitem{A2} {\sc M.A.M. Alwash}, {\it The composition conjecture for Abel equation},
    Expo. Math. {\bf 27} (2009), no. 3, 241--250.

\bibitem{AL} {\sc M.A.M. Alwash, N.G. Lloyd}, {\it Non-autonomous equations related to polynomial two-dimensional
systems},  Proc. Roy. Soc. Edinburgh Sect. A {\bf 105} (1986), 129--152.

\bibitem{BlY} {\sc M. Blinov, Y. Yomdin}, {\it Generalized center conditions and multiplicities for polinomial Abel equations of small degrees}, Nonlinearity {\bf 12} (1999), 1013--1028.

\bibitem{BFY1} {\sc M. Briskin, J.P. Fran\c{c}oise, Y. Yomdin},
{\it Center conditions, compositions of polynomials and moments on algebraic curves},  Ergodic Theory Dynam. Systems {\bf 19} (1999), 1201--1220.

\bibitem{BFY2} {\sc M. Briskin, J.P. Fran\c{c}oise, Y. Yomdin},
{\it Center conditions. II. Parametric and model center problems}, Israel J. Math. {\bf 118} (2000), 61--82.

\bibitem{BFY3} {\sc M. Briskin, J.P. Fran\c{c}oise, Y. Yomdin},
{\it Center conditions. III. Parametric and model center problems}, Israel J. Math. {\bf 118} (2000), 83--108.

\bibitem{BRY} {\sc M. Briskin, N. Roytvarf, Y. Yomdin}, {\it Center conditions at infinity for Abel differential equations}, Ann. of Math. (2) {\bf 172} (2010), no. 1, 437--483.

\bibitem{Bru} {\sc A. Brudnyi}, {\it An explicit expression for the first return map in the center
problem}, J. Differential Equations {\bf 206} (2004), no. 2, 306--314.

\bibitem{Br2} {\sc A. Brudnyi}, {\it An algebraic model for the center problem},
Bull. Sci. Math. {\bf 128} (2004), no. 10, 839--857.

\bibitem{Br} {\sc A. Brudnyi}, {\it On the center problem for ordinary differential equations},
Amer. J. Math. {\bf 128} (2006), 419--451.

\bibitem{CGM} {\sc A. Cima, A. Gasull, F. Ma\~{n}osas}, {\it Centers for trigonometric Abel equations},
 Qual. Theory Dyn. Syst. {\bf 11} (2012), no. 1, 19–37.

\bibitem{CGM2} {\sc A. Cima, A. Gasull, F. Ma\~{n}osas}, {\it A simple solution of some conjectures for Abel equations}, J. Math. Anal. Appl. {\bf 398} (2013), 477–486.

\bibitem{CGM3} {\sc A. Cima, A. Gasull, F. Ma\~{n}osas}, {\it  An explicit bound of the number of vanishing double moments forcing composition}, J. Differential Equations {\bf 255} (2013), no. 3, 339--350.

\bibitem{GGL} {\sc J. Gin\'e, M. Grau, J. Llibre}, {\it Universal centers and composition conditions},
Proc. Lond. Math. Soc. (3) {\bf 106} (2013), 481--507.

\bibitem{GGX0} {\sc J. Gin\'e, M. Grau, X. Santallusia}, {\it Composition conditions in the trigonometric Abel equations}, J. Appl. Anal. Comput. {\bf 3} (2013), no. 2, 133--144.

\bibitem{GGX} {\sc J. Gin\'e, M. Grau, X. Santallusia}, {\it Universal centers in the cubic trigonometric Abel equation}, Electron. J. Qual. Theory Differ. Equ. {\bf 2014}, No. 1, 1--7.

\bibitem{GGX2} {\sc J. Gin\'e, M. Grau, X. Santallusia}, {\it The center problem and composition condition for Abel differential equations}, Expositiones Mathematicae {\bf 34} (2016), no. 2, 210--222.
    
\bibitem{GS} {\sc J. Gin\'e, X. Santallusia}, {\it Abel differential equations admitting a certain first integral}, J. Math. Anal. Appl. 370 (2010), no. 1, 187--199.

\end{thebibliography}
\end{document}